\documentclass[10pt]{amsart}
\usepackage{amsmath, amsfonts, amssymb, amstext, amscd, amsthm}
\numberwithin{equation}{section}
\allowdisplaybreaks

\newtheorem{descent}{Theorem}[section]
\newtheorem{4-isogeny}[descent]{Proposition}
\newtheorem{mordell-weil}[descent]{Proposition}
\newtheorem{weil-pairing}[descent]{Proposition}
\newtheorem{weil-chatelet}[descent]{Proposition}
\newtheorem{selmer-sha}[descent]{Proposition}
\newtheorem{4-descent}[descent]{Theorem}
\newtheorem{2-isogeny}[descent]{Proposition}
\newtheorem{2-descent}[descent]{Proposition}

\newcommand{\sha}{\makebox[0pt][l]{\hspace{.005in}\rule{.151in}{.003in}}{\rm III}}

\begin{document}


\title{Explicit Descent via 4-Isogeny on an Elliptic Curve}

\author{Edray Herber Goins}
\address{Department of Mathematics \\ California Institute of Technology \\ Pasadena, CA 91125}
\email{goins@caltech.edu}


\begin{abstract}
We work out the complete descent via 4-isogeny for a family of rational elliptic curves with a rational point of order 4; such a family is of the form $y^2 + x \, y + a \, y = x^3 + a \, x^2$ where $\sqrt{-a} \in \mathbb Q^\times$.  In the process we exhibit the 4-isogeny and the isogenous curve, explicitly present the principal homogeneous spaces, and discuss examples by computing the rank.
\end{abstract}

\maketitle


\section{Introduction}

For each $t \in \mathbb Q^\times$ we have the rational elliptic curve
\begin{equation} E_t: \quad v^2 = u^3 + \left( t^2 + 2 \right) u^2 + u \quad \text{where} \quad [4] \left( -1, t \right) = \mathcal O. \end{equation}

\noindent This curve is studied in detail in \cite{CG1}, where, among other things, it is shown that the torsion subgroup $\mathbb T$ of $E_t(\mathbb Q)$ is completely classified as being either
\begin{equation} \mathbb T \simeq \begin{cases} \mathbb Z / 4 \mathbb Z & \text{or} \\ \mathbb Z / 2 \mathbb Z \times \mathbb Z / 4 \mathbb Z & \text{if $t = (s^2 - 1)/s$ for $s \in \mathbb Q^\times$, or} \\ \mathbb Z / 2 \mathbb Z \times \mathbb Z / 8 \mathbb Z & \text{if $t = (s^2 - 1)/s$ and $s = (r^2 - 1)/(2 \, r)$ for $r \in \mathbb Q^\times$.} \end{cases} \end{equation}

\noindent On the other hand, the rank of $E_t(\mathbb Q)$ may be computed by a descent via two-isogeny coming from the rational point $[2] \, (-1,t) = (0,0)$ of order 2, as implemented in computer packages such as \texttt{MAGMA} \cite{MR1484478}.   Unfortunately, there may be a global obstruction in computing the rank due to nontrivial elements in the Shafarevich-Tate group, so one can only find a lower bound on the rank in general.  In this exposition, I show how to compute the rank of $E_t(\mathbb Q)$ by exploiting instead the rational point $(-1,t)$ of order 4 i.e. I work out a descent via four-isogeny.

The main result of this paper is the following:
\begin{descent} Fix $t \in \mathbb Q^\times$, and denote the elliptic curves
\begin{equation} \begin{aligned} E_t: & \quad & v^2 & = u^3 + \left( t^2 + 2 \right) u^2 + u, \\ E_t': & \quad & V^2 & = U^3 - 2 \left( t^2 - 4 \right) \, U^2 + \left( t^2 + 4 \right)^2 U, \\ E_t'': & \quad  & {\mathfrak v}^2 & = {\mathfrak u}^3 + \left( t^2 - 4 \right) {\mathfrak u}^2 - 4 \, t^2 \, \mathfrak u. \end{aligned} \end{equation}
\begin{enumerate}
\item There are isogenies $\phi: E_t \to E_t'$ of degree 4 with kernel generated by $(-1,t)$, and $\varphi: E_t \to E_t''$ of degree 2 generated by $[2] \, (-1,t) = (0,0)$.
\item Upon identifying the Selmer groups as subgroups of $\mathbb Q^\times / (\mathbb Q^\times)^4$ and $\mathbb Q^\times / (\mathbb Q^\times)^2$, respectively; and Shafarevich-Tate groups as a collection of homogeneous spaces, there is an exact diagram
\begin{equation} \begin{CD} 0 @>>> \dfrac {E_t(\mathbb Q)}{\hat \phi \left( E_t'(\mathbb Q)\right)} @>{\delta'}>> S^{(\hat \phi)}(E_t' / \mathbb Q) @>>> \sha(E_t' / \mathbb Q)[\hat \phi] @>>> 0 \\ @. @VVV @VVV @VV{\eta}V \\ 0 @>>> \dfrac {E_t(\mathbb Q)}{\hat \varphi \left( E_t''(\mathbb Q)\right)} @>{\delta''}>> S^{(\hat \varphi)}(E_t'' / \mathbb Q) @>>> \sha(E_t'' / \mathbb Q)[\hat \varphi] @>>> 0 \end{CD} \end{equation} 
\noindent where we have the connecting homomorphisms
\begin{equation} \delta': \begin{aligned} (-1,t) & \mapsto - 4 \, t^2 \\ (u,v) & \mapsto u^2 + 2 \, (t^2+1) \, u + 1 - 2 \,t \, v \end{aligned} \qquad \delta'': \quad \begin{aligned} (0,0) & \mapsto  t^2 + 4  \\ (u,v) & \mapsto u \end{aligned} \end{equation}
\noindent while the maps into the Shafarevich-Tate groups send $d \mapsto \{ C_{d,t}' / \mathbb Q \}$ and $d \mapsto \{ C_{d,t}'' / \mathbb Q \}$ in terms of
\begin{equation} \begin{aligned} C_{d,t}': & \quad & d \left( w - \frac {z^2}{4 \, t^2} \right) z^2 = \left( w^2 - d \right)^2, \\ C_{d,t}'': & \quad & d \, W^2 = d^2 + \left( t^2 + 2 \right) d \, Z^2 + Z^4. \end{aligned} \end{equation}
\item The composition
\begin{equation} \begin{CD} \psi': \quad C_{d,t}' @>{\eta}>> C_{d,t}'' @>{\psi''}>> E_t \end{CD} \end{equation} 
\noindent in terms of the maps $\eta$ and $\psi''$ sending
\begin{equation} (z,w) \mapsto \left( -\frac {d \, z^2}{2 \, t \, (w^2 - d)}, \ \frac {d \, z^2 \, (w^2 + d)}{2 \, (w^2 - d)^2} \right), \qquad (Z,W) \mapsto \left( \frac {d}{Z^2}, \ - \frac {d \, W}{Z^3} \right); \end{equation}
\noindent satisfy the congruence $\delta' \bigl(  \psi'(P) \bigr) \equiv d \mod{(\mathbb Q^\times)^4}$.
\end{enumerate}
\end{descent}

I mention that the article \cite{MR97j:11027} outlines a similar algorithm using a 4-descent by considering a diagram such as
\begin{equation} \begin{CD} 0 @>>> \dfrac {E(\mathbb Q)}{4 \, E(\mathbb Q)} @>>> S^{(4)}(E / \mathbb Q) @>>> \sha(E / \mathbb Q)[4] @>>> 0 \\ @. @VVV @VVV @VVV \\ 0 @>>> \dfrac {E(\mathbb Q)}{2 \, E(\mathbb Q)} @>>> S^{(2)}(E / \mathbb Q) @>>> \sha(E / \mathbb Q)[2] @>>> 0 \end{CD} \end{equation} 

\noindent (I have learned that Tom Womack is working on implementation of a similar algorithm for a future version of \texttt{MAGMA}.)  The paper at hand may be considered a refinement of this algorithm, even though I present the algorithm for a limited family of elliptic curves possessing a 4-isogeny.

I am motivated to implement such an algorithm by the search for elliptic curves $E_t$ with $\mathbb T \simeq \mathbb Z / 2 \mathbb Z \times \mathbb Z / 8 \mathbb Z$ having ``large'' rank.  In \cite{CG1}, there is an elliptic curve with rank 3 corresponding to $r = 15/76$.  Moreover, the author, along with Garikai Campbell, conjectured that if
\begin{equation} r \in \left \{ \frac {15}{56}, \, \frac {24}{65}, \, \frac {11}{69}, \, \frac 7{88}, \, \frac {12}{97} \right \} \end{equation}

\noindent then $E_t(\mathbb Q)$ will have rank at least 4.  Randall Rathbun \cite{rathbun:01} has informed me that if we assume the validity of the conjectures of Birch and Swinnerton-Dyer then our conjecture is false; in fact the curves have rank not greater than 2.  This is apparently because the Shafarevich-Tate groups in the sequence
\begin{equation} \begin{CD} 0 @>>> \sha(E_t/\mathbb Q)[\phi] @>>> \sha(E_t /\mathbb Q)[4] @>{\phi}>> \sha(E_t'\mathbb Q)[\hat \phi] @>>> 0 \end{CD} \end{equation}

\noindent are quite large.  Explicitly, for these five values of $r$ we have $\# \sha(E_t'\mathbb Q)[\hat \phi] \geq 2^4$.

The authors would like to thank Garikai Campbell, John Cremona, Lloyd Kilford, William Stein, and Randall Rathbun for helpful conversations.

\tableofcontents


\section{Descent via Four-Isogeny}

For this section, we fix $t \in \mathbb Q^\times$ and drop the notation involving $t$. 

\subsection{The 4-Isogeny}

We begin by explicitly describing the 4-isogeny and the isogenous elliptic curve. 

\begin{4-isogeny} \label{4-isogeny} The elliptic curve $E_t$ is rationally equivalent to
\begin{equation} E: \quad y^2 + x \, y + a \, y = x^3 + a \, x^2 \quad \text{where} \quad a = -\frac 1{4 \, t^2}, \quad [4] \, (0,0) = \mathcal O; \end{equation}
\noindent and is 4-isogenous to the elliptic curve
\begin{equation} E': \quad Y^2 + X \, Y + A \, Y = X^3 + A \, X^2 \quad \text{where} \quad A = \frac {t^2 + 4}{64}. \end{equation}

\noindent Moreover, the isogeny $\phi: E \to E'$ sends $(0,0) \mapsto \mathcal O$.
\end{4-isogeny}

One can easily show that any rational elliptic curve with a rational point of order 4 is birationally equivalent to a curve in the form above for some $a \in \mathbb Q^\times$.  Indeed, see \cite[Exercise 8.13(a)]{MR87g:11070} as well as \cite[pg. 217]{MR55:7910}.  The condition $\sqrt{-a} \in \mathbb Q^\times$ is special to simplify the formulas to follow.

\begin{proof} Starting with the curve $v^2 = u^3 + \left( t^2 + 2 \right) u^2 + u$, make the substitution
\begin{equation} x = \frac {u + 1}{4 \, t^2} \quad \text{and} \quad y = \frac {v - t \, u}{8 \, t^3} \implies (-1, \, t) \mapsto (0,0). \end{equation}

\noindent Then $E$ has the Weierstrass equation above.  Following the ideas in \cite{MR45:3414} we produce an isogeny $\phi: E \to E'$ such that the cyclic subgroup 
\begin{equation} E[\phi] = \left \{ \left( 0, 0 \right), \left( -a, 0 \right), \left( 0, -a \right), \mathcal O \right \} \end{equation}

\noindent is the kernel.  We choose such a map that sends $(x,y) \mapsto (X,Y)$ in terms of
\begin{equation} \begin{aligned} X & = -A + \left( \frac {t \, (x+2 \, a) \, (2 \, y + x + a)}{8 \, x \, (x+a)} \right)^2, \\ Y & = X^2 - \left( x+a \right)^2 \left( \frac {(2 \, y + x + a) - 2 \, t \, x \, (x + 2 \, a)}{8 \, x \, (x+a)} \right)^4. \end{aligned} \end{equation}

\noindent The new curve $E'$ has a rational point of order 4, but we choose the dual isogeny $\hat \phi: E' \to E$ such that
\begin{equation} E'[\hat \phi] = \left \{ \left( -2 \, A, \ A \, \frac {2 + i \, t}4 \right), \left( - A, \ 0 \right), \ \left( -2 \, A, \ A \, \frac {2 - i \, t}4 \right), \ \mathcal O \right \} \end{equation}

\noindent where $i = \sqrt{-1}$.  This map sends $(X,Y) \mapsto (x,y)$ in terms of
\begin{equation} \label{dual-isogeny} \begin{aligned} x & = -a + \left( \frac {X \, (2 \, Y + X + A)}{2 \, t \, (X+A) \, (X+2 \, A)} \right)^2, \\ y & = x^2 - \left( \frac {2 \, (X + 2 \, A) \, (2 \, Y + X + A) - t \, X \, (X+A)}{4 \, t \, (X+A) \, (X+2 \, A)} \right)^4. \end{aligned} \end{equation}

\noindent One checks that $\hat \phi \circ \phi = [4]$ is the ``multiplication-by-4'' map on $E$, while $\phi \circ \hat \phi = [4]$ is the ``multiplication-by-4'' map on $E'$.   \end{proof}

\subsection{The Isogenous Curve}

The following proposition shows how the size of $E(\mathbb Q) / 4 \, E(\mathbb Q)$ is related to the size of $E(\mathbb Q) / \hat \phi \left( E'(\mathbb Q) \right)$.

\begin{mordell-weil} \label{mordell-weil}  Let $\phi: E \to E'$ be the isogeny in Proposition \ref{4-isogeny} and $\hat \phi: E' \to E$ denote its dual.
\begin{enumerate}
\item The following sequence is exact if and only if $\sqrt{t^2 + 4} \in \mathbb Q^\times$:
\begin{equation} \begin{CD} 0 @>>> \dfrac {E'(\mathbb Q)}{\phi \left( E(\mathbb Q) \right)} @>{\hat \phi}>> \dfrac {E(\mathbb Q)}{4 \, E(\mathbb Q)} @>>> \dfrac {E(\mathbb Q)}{\hat \phi \left( E'(\mathbb Q) \right)} @>>> 0. \end{CD} \end{equation}
\item Both $E(\mathbb Q)$ and $E'(\mathbb Q)$ have the same rank.
\item The torsion subgroup $\mathbb T'$ of $E'(\mathbb Q)$ is $\mathbb Z / 8 \, \mathbb Z$ when $\sqrt{t + \sqrt{t^2 + 4}} \in \mathbb Q^\times$, and $\mathbb Z / 4 \, \mathbb Z$ otherwise.
\end{enumerate}
\end{mordell-weil}

\begin{proof} The dual isogeny induces the exact sequence
\begin{equation} 0 \to \frac {E'(\mathbb Q)[\hat \phi]}{\phi \left( E(\mathbb Q)[4] \right)} \to \frac {E'(\mathbb Q)}{\phi \left( E(\mathbb Q) \right)} \to \frac {E(\mathbb Q)}{4 \, E(\mathbb Q)} \to \frac {E(\mathbb Q)}{\hat \phi \left( E'(\mathbb Q) \right)} \to 0. \end{equation}

\noindent We compute the first term in this sequence.  From above, $E'(\mathbb Q)[\hat \phi] = \left \{ (-A,0), \, \mathcal O \right \}$ has order 2.  $E[4]$ is generated by $P = (0,0)$ and
\begin{equation} Q = \left( -\frac {2 \, a}{1 + \alpha}, \ \frac {8 \, a^3}{(1 + \alpha)^2 (1 + i \, \alpha)} \right) \quad \text{where} \quad \alpha = \sqrt[4]{ \frac {t^2 + 4}{t^2} }. \end{equation}

\noindent One verifies that $\phi(P) = \mathcal O$ while $\phi([2]Q) = (-A, 0)$.  Clearly $Q \not \in E(\mathbb Q)[4]$, while
\begin{equation} [2] \, Q = \left( - \frac {1 + \alpha^2}8, \ \frac {(1 + \alpha^2)^2}{32} \right) \in E(\mathbb Q)[4] \iff \sqrt{t^2 + 4} \in \mathbb Q^\times. \end{equation}  

\noindent Hence $E'(\mathbb Q)[\hat \phi] = \phi \left( E(\mathbb Q)[4] \right)$ if and only if $t^2 + 4$ is a square.

Second, the statement that isogenous elliptic curves have the same rank is well-known: As $E'(\mathbb Q) / \phi \left( E(\mathbb Q) \right)$ is finite the rank of $E'(\mathbb Q)$ cannot be greater than the rank of $E(\mathbb Q)$.  Similarly, $E(\mathbb Q) / \hat \phi \left( E'(\mathbb Q) \right)$ is finite so we must have equality with the rank.

As for the torsion, it is easy to see that $(0,0) \in E'(\mathbb Q)$ is a rational point of order 4.  By Mazur's Theorem, either $E'[2] \subseteq E'(\mathbb Q)$ or else the torsion subgroup is of the form $\mathbb Z /  4\, n \, \mathbb Z$ for $n = 1, \, 2, \, 3$.  We begin by showing that not all of the 2-torsion is rational. The points $(X,Y)$ of order 2 satisfy $2 \, Y + X + A = 0$, so that
\begin{equation} \left( X + A \right) \left( 4 \, X^2 + X + A \right) = -4 \left( Y^2 + X \, Y + A \, Y - X^3 - A \, X^2 \right) = 0. \end{equation}

\noindent $(-A,0)$ is one rational point of order 2 corresponding to the linear factor, but the quadratic factor has discriminant $- 1$ so that polynomial has no rational roots.

Next we show that there are no rational points of order 3.  Assume that $(X,Y)$ is such a point.  One computes that the 3-division polynomial of $E'$ is
\begin{equation} \psi_3(X,t) = 3 \, X^4 + (1 + 4 \, A) \, X^3 + 3 \, A \, X^2 + 3 \, A^2 \, X + A^3. \end{equation}

\noindent  We make the substitution
\begin{equation} \sigma = - \frac {3 \, X + 4 \, A}X \quad \text{and} \quad \tau = \frac {3 \, X + 4 \, A}{X + 2 \, A} \, t. \end{equation}

\noindent Note that $X \, (X + 2 \, A) \neq 0$ because otherwise $(X,Y)$ would be a point of order 4.  Also, $3 \, X + 4 \, A \neq 0$: if not, then $(X,Y) = (-4 \, A/3, \ -A \, \beta / 9)$ in terms of a root of the quadratic $\beta^2 + 3 \, \beta + 48 \, A$.  However, this quadratic has discriminant $-3 \, (t^2 + 1)$ which is not a square, so $(X,Y)$ is not a rational point.  Hence $\sigma$ and $\tau$ are well-defined and nonzero.  Solving for $X$ and $t$ transforms the 3-division polynomial as 
\begin{equation} \psi_3 \left(  - \frac {4 \, A}{\sigma + 3}, \ \frac {\tau \, (\sigma + 1)}{2 \, \sigma} \right)  = \frac {A^3 \, (\sigma+1)^2}{\sigma \, (\sigma+3)^4} \cdot \left[ \tau^2 - \sigma \, (\sigma + 1) \, (\sigma - 3) \right] \end{equation}

\noindent so that $(\sigma, \, \tau)$ is a rational point on the elliptic curve $\tau^2 = \sigma \, (\sigma + 1) \, (\sigma - 3)$.  However, one checks that the only such (affine) rational points satisfy $\tau = 0$, a contradiction.  

Finally, we are in the case where we distinguish between $\mathbb Z / 4 \, \mathbb Z$ and $\mathbb Z / 8 \, \mathbb Z$.  Say that $R = (X,Y)$ is a rational point of order 8.  Then $[2] \, R$ is a point of order 4, so without loss of generality say $[2] \, R = P$.  By considering the $X$-coordinates we have
\begin{equation} X \left( [2] \, R \right) = \frac {X^4 - A \, (X+A)}{(X+A) \, (4 \, X^2 + X + A)} = 0 \implies X^2 - \sqrt{A} \, ( X + A) = 0. \end{equation}

\noindent Hence $\sqrt{A} \in \mathbb Q^\times$.  In fact, make the substitution
\begin{equation} \gamma = \sqrt{t + \sqrt{t^2 + 4}} \implies t = \frac {\gamma^4 - 4}{2 \, \gamma^2} \quad \text{and} \quad A = \left( \frac {\gamma^4 + 4}{16 \, \gamma^2} \right)^2. \end{equation}

\noindent The discriminant of the quadratic $X^2 - \sqrt{A} \, ( X + A)$, namely $\gamma^2$, must be a square since its roots are rational; hence $\gamma \in \mathbb Q^\times$.  Then the rational point of order 8 is 
\begin{equation} R = \left( \frac {(\gamma^4 + 4) \, (\gamma^2 + 2 \, \gamma + 2)}{64 \, \gamma^3}, \ \frac {(\gamma^4 + 4)^2 \, (\gamma^2 + 2 \, \gamma + 2)}{1024 \, \gamma^5} \right). \end{equation}

\noindent The converse is clear i.e. if $[8]R = \mathcal O$ then $\sqrt{t + \sqrt{t^2 + 4}} \in \mathbb Q^\times$.  \end{proof}

\subsection{The Connecting Homomorphism}

There is a natural bijection between the crossed homomorphisms $H^1 \left( \text{Gal}(\overline {\mathbb Q}/\mathbb Q), \, E \right)$ and $WC \left( E/\mathbb Q \right)$, the collection of equivalence classes of homogeneous spaces for $E$, usually called the Weil-Ch{\^a}telet group.  We will find it more convenient to work with the isogenous curve $E'$:  The exact sequence
\begin{equation} \begin{CD} 0 @>>> E'[\hat \phi] @>>> E' @>{\hat \phi}>> E @>>> 0 \end{CD} \end{equation}

\noindent implies, through Galois cohomology, the exact sequence
\begin{equation} \begin{CD} 0 @>>> \dfrac {E(\mathbb Q)}{\hat \phi \left( E'(\mathbb Q)\right)} @>{\delta_E}>> H^1 \left( \text{Gal}(\overline {\mathbb Q}/\mathbb Q), \, E'[\hat \phi] \right) @>>> WC \left( E'/\mathbb Q \right) \end{CD} \end{equation}

\noindent via the connecting homomorphism $\delta_E$.  We make such a homomorphism explicit by considering a generalization of the Weil pairing.

\begin{weil-pairing} \label{weil-pairing}
\begin{enumerate}
\item There is a pairing $e_{\phi}: E'[\hat \phi] \times E[\phi] \to \mu_4$ that is bilinear, alternating, non-degenerate, and Galois invariant; which satisfies
\begin{equation} e_{\phi}(T',T) = \sqrt{-1} \quad \text{when} \quad T' = \left( -2 \, A, \ A \, \frac {2 + \sqrt{-1} \, t}4 \right), \quad T = (0,0). \end{equation}
\item The pairing $b: E(\mathbb Q) / \hat \phi \left( E'(\mathbb Q)\right) \times E[\phi] \to \mathbb Q^\times / (\mathbb Q^\times)^4$ defined by
\begin{equation} \biggl( (x,y), \ [m] \, (0,0) \biggr) \mapsto \begin{cases} a^{-m} & (x,y) = (0,0), \\ \left( x^2 - y \right)^m & \text{otherwise;} \end{cases} \end{equation}
\noindent is bilinear and non-degenerate on the left.
\item Upon identifying $E'[\hat \phi] \simeq \mu_4$ as $\text{Gal}(\overline {\mathbb Q}/\mathbb Q)$-modules we have
\begin{equation} H^1 \left( \text{Gal}(\overline {\mathbb Q}/\mathbb Q), \, E'[\hat \phi] \right) \simeq \mathbb Q^\times / (\mathbb Q^\times)^4. \end{equation}
\item The composition of these maps
\begin{equation} \begin{CD} \dfrac {E(\mathbb Q)}{\hat \phi \left( E'(\mathbb Q)\right)} \times E[\phi] @>{b}>> \dfrac {\mathbb Q^\times}{(\mathbb Q^\times)^4} @>{\delta_{\mathbb Q}}>> H^1 \left( \text{Gal}(\overline {\mathbb Q}/\mathbb Q), \, E'[\hat \phi] \right) \end{CD}  \end{equation}
\noindent relates the connecting homomorphisms by
\begin{equation} e_{\phi} \bigl( \delta_E(P)(\sigma), \ T \bigr) = \delta_{\mathbb Q} \bigl( b(P,T) \bigr) (\sigma) \quad \text{for all $\sigma \in \text{Gal}( \overline {\mathbb Q} / \mathbb Q)$.} \end{equation}
\noindent In particular, we may identify $\delta_E: P \mapsto \delta_{\mathbb Q} \bigl( b(P,T) \bigr)$ when $T = (0,0)$.
\end{enumerate} \end{weil-pairing}

\begin{proof} We follow the construction of the Weil pairing in \cite[\S 8, pgs. 95 -- 99]{MR87g:11070}, but we make the formulas explicit.  Consider the function $f \in \mathbb Q(E)$ defined as $f(x,y) = x^2 - y$.  One checks that
\begin{equation} x^2 - y = 0 \implies x^4 = y^2 + x \, y + a \, y - x^3 - a \, x^2 = 0 \implies (x,y) = (0,0); \end{equation}

\noindent so that the divisor of this function is $\text{div}(f) = 4 \left( (0,0) \right) - 4(\mathcal O)$.  Similarly, we consider the function $g \in \mathbb Q(E')$ defined by
\begin{equation} g(X,Y) = \frac {2 \, (X + 2 \, A) \, (2 \, Y + X + A) - t \, X \, (X+A)}{4 \, t \, (X+A) \, (X+2 \, A)}; \end{equation}

\noindent the formulas in \eqref{dual-isogeny} show that $f \circ \hat \phi = g^4$.  For $P' \in E'(\overline {\mathbb Q})$ and $T' \in E'[\hat \phi]$ we have
\begin{equation} g \left( P' \oplus T' \right)^4 =f \left( \hat \phi(P') \oplus \hat \phi(T') \right) = f \left( \hat \phi(P') \oplus \mathcal O \right) = g \left( P' \right)^4. \end{equation}

\noindent We define the pairing $e_{\phi} : E'[\hat \phi] \times E[\phi] \to \mu_4$ by
\begin{equation} e_{\phi} \left( T', \, T \right) = \left( \frac {g \left( P' \oplus T' \right)}{g \left( P' \right)} \right)^m \quad \text{where} \quad T = [m] \, (0,0). \end{equation}

\noindent The bilinear, alternating, non-degenerate, and Galois invariant properties follow from the arguments given in  \cite[\S 8, pgs. 95 -- 99]{MR87g:11070}.  In particular, one verifies that
\begin{equation} T' = \left( -2 \, A, \ A \, \frac {2 + \sqrt{-1} \, t}4 \right), \quad T = (0,0) \implies e_{\phi} \left(T', \, T \right) = \frac {g(P' \oplus T')}{g(P')} = \sqrt{-1}. \end{equation}

Consider the map $b: E(\mathbb Q) \times E[\phi] \to \mathbb Q^\times / (\mathbb Q^\times)^4$.  First we show it is well-defined.  As shown above, $x^2 - y = 0$ if and only if $(x,y) = (0,0)$.  On the other hand, we have by assumed bilinearity
\begin{equation} b \biggl( (0,0), \ [m] \, (0,0) \biggr) = b \biggl( [-1] \, (0,0), \ [m] \, (0,0) \biggr)^{-1}  = a^{-m}. \end{equation}

\noindent Now we show that this map can be extended to the quotient.  When $(x,y) = \hat \phi (X,Y)$ then $x^2 - y$ is a fourth-power by \eqref{dual-isogeny} so that the map extends to a well-defined map $b: E(\mathbb Q) / \hat \phi \left( E'(\mathbb Q) \right) \times E[\phi] \to \mathbb Q^\times / (\mathbb Q^\times)^4$.  The fact that $b$ is a non-degenerate pairing follows from the same arguments as with the Weil pairing above.  

Fix $d \in \mathbb Q^\times$, and consider the map $\xi: \text{Gal}(\overline {\mathbb Q}/\mathbb Q) \to E'[\hat \phi]$ defined by
\begin{equation} \xi: \quad \sigma \mapsto \xi_\sigma = [m] \left( -2 \, A, \ A \, \frac {2 + i \, t}4 \right) \quad \text{whenever} \quad \frac {\sigma (\sqrt[4]{d})}{\sqrt[4]{d}} = (-1)^{m/2}. \end{equation}

\noindent The map $\xi$ is a 1-cocycle, so has a cohomology class $\{ \xi \} \in H^1 \left( \text{Gal}(\overline {\mathbb Q}/\mathbb Q), \, E'[\hat \phi] \right)$.  Identify translation by this point with the Galois action i.e. $(X,Y)^\sigma = (X,Y) \oplus \xi_\sigma$.  From the map $d \mapsto d^4$ we find the well-known Kummer sequence
\begin{equation} \begin{CD} 1 @>>> E'[\hat \phi] \simeq \mu_4 @>>> \overline {\mathbb Q}^\times @>>> \overline {\mathbb Q}^\times @>>> 1 \end{CD} \end{equation}

\noindent of $\text{Gal}(\overline {\mathbb Q}/\mathbb Q)$-modules, so from Galois cohomology we find the exact sequence
\begin{equation} \begin{CD} 1@>>> \dfrac {\mathbb Q^\times}{(\mathbb Q^\times)^4} @>{\delta_{\mathbb Q}}>> H^1 \left( \text{Gal}(\overline {\mathbb Q}/\mathbb Q), \, E'[\hat \phi] \right) @>>> H^1 \left( \text{Gal}(\overline {\mathbb Q}/\mathbb Q), \, \overline {\mathbb Q}^\times \right). \end{CD} \end{equation}

\noindent where $\delta_{\mathbb Q}: d \mapsto \{ \xi \}$. Hilbert's Theorem 90 states the last term is trivial; hence the third part of the proposition holds.

To show the fourth part of the proposition, fix $P \in E(\mathbb Q)$, choose $P' \in E'(\overline {\mathbb Q})$ satisfying $P = \hat \phi(P')$, and denote $\xi = \delta_E(P)$.  For any $\sigma \in \text{Gal}( \overline {\mathbb Q} / \mathbb Q)$ we have
\begin{equation} \begin{aligned} e_{\phi} \bigl( \delta_E(P)(\sigma), \, T \bigr) & = e_{\phi} \bigl( \xi_\sigma, \, [m] (0,0) \bigr) = \left( \frac {g \left( P' \oplus \xi_\sigma \right)}{g(P')} \right)^m = \left( \frac {g(P')^{\sigma}}{g(P')} \right)^m \\ & = \bigl( \delta_{\mathbb Q}(d) (\sigma) \bigr)^m = \delta_{\mathbb Q} \bigl( b(P,T) \bigr) (\sigma). \end{aligned} \end{equation}

\noindent where we have identified $E'[\hat \phi] \simeq \mu_4$ and set $d = f(P) \equiv b \left( P, \, (0,0) \right)$. \end{proof}

\subsection{The Weil-Ch{\^a}telet Group}

The full collection of homogeneous spaces in the Weil-Ch{\^a}telet group is too large for interest in this exposition, so instead we consider the image $WC(E'/\mathbb Q)[\hat \phi]$ of the crossed homomorphisms studied above.  Recall that $WC(E'/\mathbb Q)$ consists of equivalence classes $\{ C'/\mathbb Q \}$ of certain rational curves, where a class is trivial if and only  if $C'(\mathbb Q)$ is nonempty. 

\begin{weil-chatelet} \label{weil-chatelet}
\begin{enumerate}
\item The composite map
\begin{equation} \begin{CD} \dfrac {\mathbb Q^\times}{(\mathbb Q^\times)^4} @>{\delta_{\mathbb Q}}>> H^1 \left( \text{Gal}(\overline {\mathbb Q}/\mathbb Q), \, E'[\hat \phi] \right) @>>> WC \left( E'/\mathbb Q \right) \end{CD} \end{equation}

\noindent sends $d \mapsto \{ C_d'/\mathbb Q \}$, where
\begin{equation} C_d': \quad d \left( w + a \, z^2 \right) z^2  = \left( w^2 - d \right)^2. \end{equation}
\item There is an isomorphism $\theta: C_d' \to E'$ defined over $\overline {\mathbb Q}$ such that
\begin{equation} \begin{CD} \psi': \quad C_d' @>{\theta}>> E' @>{\hat \phi}>> E, \end{CD} \qquad (z,w) \mapsto \left( \frac {w}{z^2}, \ \frac {w^2 - d}{z^4} \right). \end{equation}
\end{enumerate}
\end{weil-chatelet}

\begin{proof}  We construct a homogeneous space corresponding the cohomology class $\{ \xi \} \in H^1 \left( \text{Gal}(\overline {\mathbb Q}/\mathbb Q), \, E'[\hat \phi] \right)$ following the exposition in \cite[\S 3, pgs. 287 -- 296]{MR87g:11070}.  Fix $d \in \mathbb Q^\times$, and consider the functions $z, \, w \in \overline {\mathbb Q}(E')$ defined implicitly by
\begin{equation}  z =  \frac {\sqrt[4]{d}}{g(P')} \quad \text{and} \quad - \frac zw = \frac {1/\sqrt[4]{d}}{g([-1]P')} \quad \text{where} \quad P' = (X,Y). \end{equation}

\noindent It is easy to check that the point $(z,w)$ is on $C_d'$ as defined above.  We show these functions are Galois invariant.  To this end, choose $\sigma \in \text{Gal}(\overline {\mathbb Q} / \mathbb Q)$, and say that $\sigma(\sqrt[4]{d}) = (-1)^{m/2} \, \sqrt[4]{d}$.  For any integer $n$ we have
\begin{equation} \sigma \left( \frac {\sqrt[4]{d^n}}{g([n] P')} \right) = \frac {(-1)^{mn/2} \, \sqrt[4]{d^n}}{g([n] P'^\sigma)} = \frac {(-1)^{mn/2} \, \sqrt[4]{d^n}}{g([n] P' \oplus T')} = \frac {(-1)^{mn/2}}{e_{\phi}(T',T)} \ \frac {\sqrt[4]{d^n}}{g([n] P')}\end{equation}

\noindent where $T' = [n] \, \xi_\sigma$ and $T = (0,0)$.  It follows from Proposition \ref{weil-pairing} that $e_{\phi}(T',T) = (-1)^{mn/2}$.  Hence the equations defining $z$ and $w$ are Galois invariant, so that $z$ and $w$ themselves are Galois invariant.

This defines a map $E' \to C_d'$ sending $(X, \, Y) \mapsto (z,w)$ in terms of
\begin{equation}  \begin{aligned} z & = \sqrt[4]{d} \ \frac {4 \, t \, (X+A) \, (X+2 \, A)}{2 \, (X + 2 \, A) \, (2 \, Y + X + A) - t \, X \, (X+A)}, \\ w & = \sqrt{d} \ \frac {2 \, (X + 2 \, A) \, (2 \, Y + X + A) + t \, X \, (X+A)}{2 \, (X + 2 \, A) \, (2 \, Y + X + A) - t \, X \, (X+A)}. \end{aligned} \end{equation}

\noindent The inverse map $\theta: C_d' \to E'$ -- which can be found by considering $g(P') \pm g([-1]P')$ -- sends $(z,w) \mapsto (X,Y)$ in terms of
\begin{equation} \label{theta} X = 4 \, A \ \frac {w - \sqrt{d}}{\sqrt[4]{d} \, z - 2 \, (w - \sqrt{d})}, \quad Y = (X + A) \ \frac {-\sqrt[4]{d} \, z + t \, (w + \sqrt{d})}{2 \, \sqrt[4]{d} \, z}. \end{equation}

\noindent One checks using the formulas in \eqref{dual-isogeny} that the composition $\psi' = \hat \phi \circ \theta$ is the map $(z, \, w) \mapsto \left( w / z^2, \ (w^2 - d) / z^4 \right)$. 
\end{proof}

\subsection{Selmer and Shafarevich-Tate Groups}

Combining the exact sequences introduced above, we have the exact diagram
\begin{equation} \label{4-selmer} \begin{CD} 0 @>>> \dfrac {E'(\mathbb Q)}{\phi \left( E(\mathbb Q)\right)} @>>> S^{(\phi)}(E / \mathbb Q) @>>> \sha(E / \mathbb Q)[\phi] @>>> 0 \\ @. @V{\hat \phi}VV @VVV @VVV \\ 0 @>>> \dfrac {E(\mathbb Q)}{4 \, E(\mathbb Q)} @>>> S^{(4)}(E / \mathbb Q) @>>> \sha(E / \mathbb Q)[4] @>>> 0 \\ @. @VVV @VVV @VV{\phi}V \\ 0 @>>> \dfrac {E(\mathbb Q)}{\hat \phi \left( E'(\mathbb Q)\right)} @>{\delta_E}>> S^{(\hat \phi)}(E' / \mathbb Q) @>>> \sha(E' / \mathbb Q)[\hat \phi] @>>> 0 \\ @. @VVV @VVV @VVV \\ @. 0 @. 0 @. 0 \end{CD} \end{equation} 

\noindent where the Selmer and Shafarevich-Tate groups are defined as
\begin{equation} \begin{aligned} S^{(\phi)}(E/\mathbb Q) & = \text{ker} \left \{ H^1 \left( \text{Gal}(\overline {\mathbb Q}/\mathbb Q), \, E[\phi] \right) \to \prod_{\nu} WC \left( E / \mathbb Q_\nu \right) \right \}, \\ \sha(E/\mathbb Q) & = \text{ker} \left \{ WC \left( E / \mathbb Q \right) \to \prod_{\nu} WC \left( E / \mathbb Q_\nu \right) \right \}; \end{aligned} \end{equation}

\noindent taking the product over all places $\nu$ of $\mathbb Q$.  One can compute the rank of $E(\mathbb Q)$ with the explicit map in Proposition \ref{weil-chatelet} -- assuming that $\sha(E'/\mathbb Q)[\hat \phi]$ is trivial.

\begin{selmer-sha} \label{selmer-sha} Let $\Sigma$ be the finite set of places of $\mathbb Q$ consisting of $2$, $\infty$, and those primes occurring in the  factorization of $t$ and $t^2 + 4$.  We have
\begin{equation} \begin{aligned} S^{(\hat \phi)}(E'/\mathbb Q) & \simeq \left \{ \left. d = \pm \prod_{p \in \Sigma} p^{a_p} \in \frac {\mathbb Q^\times}{(\mathbb Q^\times)^4} \, \right \vert \,  \text{$C_d'(\mathbb Q_\nu) \neq \emptyset$ for all $\nu \in \Sigma$} \right \}, \\ \sha(E'/\mathbb Q)[\hat \phi] & \simeq \biggl \{ \biggl. \{ C_d'/\mathbb Q \} \in WC(E'/\mathbb Q)[\hat \phi] \, \biggr \vert \,  \text{$C_d'(\mathbb Q_\nu) \neq \emptyset$ for all $\nu \in \Sigma$} \biggr \}. \end{aligned} \end{equation}
\end{selmer-sha}

\begin{proof} We may identify the Selmer group with a subgroup of $\mathbb Q / (\mathbb Q^\times)^4$ by the isomorphism in Proposition \ref{weil-pairing}.  The subgroup $WC(E'/\mathbb Q)[\hat \phi] \subseteq WC(E'/\mathbb Q)$ consists of the equivalence classes $\{ C_d'/\mathbb Q \}$.  By definition, a class $\{ C_d'/\mathbb Q_\nu \} \in WC(E'/\mathbb Q_\nu)$ is trivial if and only if $C_d'(\mathbb Q_\nu)$ is nonempty. \end{proof}

\subsection{Complete Descent via 4-Isogeny}

We collect the results proved thus far into one statement.

\begin{4-descent} \label{4-descent} Fix $t \in \mathbb Q^\times$, and let $E_t$ and $E_t'$ denote the elliptic curves
\begin{equation} E_t: \ v^2 = u^3 + \left( t^2 + 2 \right) u^2 + u, \quad E_t': \ V^2 = U^3 - 2 \left( t^2 - 4 \right) \, U^2 + \left( t^2 + 4 \right)^2 U. \end{equation}

\begin{enumerate}
\item There is an isogeny $\phi: E_t \to E_t'$ of degree 4 with kernel generated by $(-1,t)$.
\item Upon identifying the Selmer and Shafarevich-Tate groups as in Proposition \ref{selmer-sha}, there is an exact sequence
\begin{equation} \begin{CD} 0 @>>> \dfrac {E_t(\mathbb Q)}{\hat \phi \left( E_t'(\mathbb Q)\right)} @>{\delta'}>> S^{(\hat \phi)}(E_t' / \mathbb Q) @>>> \sha(E_t' / \mathbb Q)[\hat \phi] @>>> 0 \end{CD} \end{equation}

\noindent where $\delta'$ sends $[m] \, (-1,t) \mapsto \left( - 4 \, t^2 \right)^m$ and $(u,v) \mapsto u^2 + 2 \, (t^2+1) \, u + 1 - 2 \,t \, v$ otherwise; while the second map sends $d \mapsto \{ C_{d,t}' / \mathbb Q \}$ in terms of
\begin{equation} C_{d,t}': \quad d \left( w - \frac {z^2}{4 \, t^2} \right) z^2 = \left( w^2 - d \right)^2. \end{equation}

\item The map $\psi': C_{d,t}' \to E_t$ sending
\begin{equation} (z,w) \mapsto \left( \frac {4 \, t^2 \, (w^2 - d)^2}{d \, z^4}, \ \frac {4 \, t^3 \, (w^2 - d) \, (w^2 + d)}{d \, z^4} \right) \end{equation} 

\noindent satisfies the congruence $\delta' \bigl(  \psi'(P) \bigr) \equiv d \mod{(\mathbb Q^\times)^4}$.
\end{enumerate}
\end{4-descent}

\begin{proof} The curves $E_t$ and $E_t'$ are related to $E$ and $E'$ by the transformation
\begin{equation} x = \frac {u + 1}{4 \, t^2}, \quad y = \frac {v - t \, u}{8 \, t^3} \quad \text{and} \quad X = \frac {U - (t^2 + 4)}{64}, \quad Y = \frac {V - 4 \, U}{512}. \end{equation}

\noindent The first statement follows from Proposition \ref{4-isogeny}.  The image of $\delta$ in the second statement follows from Proposition \ref {weil-pairing}, Proposition \ref{weil-chatelet}, and the congruence
\begin{equation} x^2 - y  \equiv \left( u^2 + 2 \, (t^2 + 1) \, u + 1 \right) - 2 \, t \, v \mod {(\mathbb Q^\times)^4}; \end{equation}

\noindent while the image of the second map follows from Proposition \ref{weil-chatelet}.  The third statement follows from Proposition \ref{weil-chatelet} where we write the map $\psi'$ in terms of $u$ and $v$ rather than $x$ and $y$.  From the identity
\begin{equation} (u,v) = \psi'(z,w) \implies \left( u^2 + 2 \, (t^2 + 1) \, u + 1 \right) - 2 \, t \, v = d \, \left( \frac {2 \, t}z \right)^4 \end{equation}

\noindent we find the congruence $\left( u^2 + 2 \, (t^2 + 1) \, u + 1 \right) - 2 \, t \, v \equiv d \mod{(\mathbb Q^\times)^4}$. \end{proof}


\section{Descent via Two-Isogeny}

The descent algorithm associated to the rational point of order 4 described in the previous section can be refined to exploit the 2-isogeny coming from the doubling of this point.  

\subsection{The 2-Isogenies}

We begin by factoring the 4-isogeny.

\begin{2-isogeny} \label{2-isogeny}  Let $\phi : E \to E'$ as in Proposition \ref{4-isogeny}, and denote the curve
\begin{equation} E'' : \qquad {\mathfrak v}^2 = {\mathfrak u} \left( \mathfrak u - 4 \right) \left( \mathfrak u + t^2 \right). \end{equation}

\noindent There exist rational 2-isogenies $\varphi: E \to E''$ and $\eta: E' \to E''$ such that $\phi = \hat \eta \circ \varphi$.
\end{2-isogeny}

\begin{proof}  The isogeny $\phi$ and the ``multiplication-by-2'' map on $E$ induce the following exact diagram:
\begin{equation} \begin{CD} @. 0 @. 0 @. 0 \\ @. @VVV @VVV @VVV \\ 0 @>>> E[\phi] \cap E[2] @>>> E[2] @>{\phi}>> \phi(E[2]) @>>> 0 \\ @. @VVV @VVV @VVV \\ 0 @>>> E[\phi] @>>> E @>{\phi}>> E' \simeq E / E[\phi] @>>> 0 \\ @. @VV{[2]}V @VV{[2]}V @VVV \\ 0 @>>> 2 \, E[\phi] @>>> E @>>> E'' := E / 2 \, E[\phi] @>>> 0 \\ @. @VVV @VVV @VVV \\ @. 0 @. 0 @. 0 \end{CD} \end{equation}

\noindent We will construct isogenies $\varphi: E \to E''$ and $\eta: E' \to E''$ with kernels
\begin{equation} E[\varphi] = 2 \, E[\phi] = \bigl \{ (-a,0), \ \mathcal O \bigr \} \quad \text{and} \quad E'[\eta] = \phi(E[2]) = \bigl \{ (-A,0), \ \mathcal O \bigr \}. \end{equation}

\noindent To this end, assume without loss of generality that $E = E_t$ and $E' = E_t'$ as in Theorem \ref{4-descent}, where now $(u,v) = (0,0)$ and $(U,V) = (0,0)$ are the rational points of order 2.  Following the exposition in \cite{MR87g:11070}, we define the isogenies
\begin{equation} \begin{aligned} \varphi : & \quad E \to E'', & \qquad (u,v) & \mapsto \left( \frac {(u+1)^2}{u},  \ \frac {1 - u^2}{u^2} \, v \right); \\ \eta: & \quad E' \to E'', & \qquad (U,V) & \mapsto \left( \frac {V^2}{4 \, U^2},  \ \frac {(t^2+4)^2 - U^2}{8 \, U^2} \, V \right); \end{aligned} \end{equation}

\noindent with dual isogenies
\begin{equation} \begin{aligned} \hat \varphi: & \quad E'' \to E, & \qquad (\mathfrak u, \mathfrak v) & \mapsto \left( \frac {{\mathfrak v}^2}{4 \, (\mathfrak u + t^2)^2}, \ - \frac {-4 \, t^2 + 2 \, t^2 \, \mathfrak u + {\mathfrak u}^2}{8 \, (\mathfrak u + t^2)^2} \, \mathfrak v \right); \\ \hat \eta : & \quad E'' \to E', & \qquad (\mathfrak u, \mathfrak v) & \mapsto \left( \frac {{\mathfrak v}^2}{{\mathfrak u}^2}, \ \frac {-4 \, t^2 - {\mathfrak u}^2}{{\mathfrak u}^2} \, \mathfrak v \right). \end{aligned} \end{equation}

\noindent One checks that $\hat \varphi \circ \varphi = [2]$ and $\hat \eta \circ \eta = [2]$ are the ``multiplication-by-2'' maps on $E$, and $E'$, respectively. Moreover, one checks that $\phi = \hat \eta \circ \varphi$ and $\hat \phi = \hat \varphi \circ \eta$.  \end{proof}

\subsection{Complete Descent via 2-Isogeny}

The previous proposition implies the following exact diagram as a refinement of the diagram in \eqref{4-selmer}:
\begin{equation} \label{descent-diagram} \begin{CD} @. 0 @. 0 @. 0 \\ @. @VVV @VVV @VVV \\ 0 @>>> \dfrac {E''(\mathbb Q)[\hat \varphi]}{\eta ( E'(\mathbb Q)[\hat \phi] )} @>{\sim}>> \left \{ \pm 1 \right \} @>>> 0 @>>> 0 \\ @. @V{\hat \varphi}VV @VVV @VVV \\ 0 @>>> \dfrac {E''(\mathbb Q)}{\eta \left( E'(\mathbb Q)\right)} @>>> S^{(\eta)}(E' / \mathbb Q) @>>> \sha(E' / \mathbb Q)[\eta] @>>> 0 \\ @. @V{\hat \varphi}VV @VVV @VVV \\ 0 @>>> \dfrac {E(\mathbb Q)}{\hat \phi \left( E'(\mathbb Q)\right)} @>{\delta_E}>> S^{(\hat \phi)}(E' / \mathbb Q) @>>> \sha(E' / \mathbb Q)[\hat \phi] @>>> 0 \\ @. @VVV @VVV @VV{\eta}V \\ 0 @>>> \dfrac {E(\mathbb Q)}{\hat \varphi \left( E''(\mathbb Q)\right)} @>>> S^{(\hat \varphi)}(E'' / \mathbb Q) @>>> \sha(E'' / \mathbb Q)[\hat \varphi] @>>> 0 \\ @. @VVV @VVV @VVV \\ @. 0 @. 0 @. 0 \end{CD} \end{equation} 

\noindent The first (nontrivial) row in the diagram above comes from the explicit relation
\begin{equation} E''(\mathbb Q)[\hat \varphi] = \left \{ (-t^2,0), \, \mathcal O \right \} \qquad \text{yet} \qquad (\mathfrak u, \mathfrak v) = \eta(U,V) \implies \mathfrak u = \frac {V^2}{4 \, U^2} > 0; \end{equation}

\noindent while the second and fourth rows are related to a descent via 2-isogeny.  We recall the main results for the latter:

\begin{2-descent} \label{2-descent} Fix $t \in \mathbb Q^\times$, denote $E_t$ as in Theorem \ref{4-descent}, and denote $E_t'' = E''$ as in Proposition \ref{2-isogeny}.

\begin{enumerate}
\item There is an exact sequence
\begin{equation} \begin{CD} 0 @>>> \dfrac {E_t(\mathbb Q)}{\hat \varphi \left( E_t''(\mathbb Q)\right)} @>{\delta''}>> S^{(\hat \varphi)}(E_t'' / \mathbb Q) @>>> \sha(E_t'' / \mathbb Q)[\hat \varphi] @>>> 0 \end{CD} \end{equation}

\noindent where $\delta''$ sends $[m] \, (0,0) \mapsto \left( t^2 + 4 \right)^m$ and $(u,v) \mapsto u$ otherwise; while the second map sends $d \mapsto \{ C_{d,t}'' / \mathbb Q \}$ in terms of
\begin{equation} C_{d,t}'': \quad d \, W^2 = d^2 + \left( t^2 +2 \right) d \, Z^2 + Z^4.  \end{equation}
\item The map $\psi': C_{d,t}' \to E_t$ in Theorem \ref{4-descent} factors as the composition
\begin{equation} \begin{CD} C_{d,t}' @>{\eta_\ast}>> C_{d,t}'' @>{\psi''}>> E_t \end{CD} \end{equation} 
\noindent where $\eta_\ast$ and $\psi''$ are the maps sending
\begin{equation} (z,w) \mapsto \left( -\frac {d \, z^2}{2 \, t \, (w^2 - d)}, \ \frac {d \, z^2 \, (w^2 + d)}{2 \, (w^2 - d)^2} \right), \qquad (Z,W) \mapsto \left( \frac {d}{Z^2}, \ - \frac {d \, W}{Z^3} \right). \end{equation}
\end{enumerate}
\end{2-descent}

In the Proposition above we identify, as in Proposition \ref{selmer-sha}, the Selmer and Shafarevich-Tate groups as
\begin{equation} \begin{aligned} S^{(\hat \varphi)}(E''/\mathbb Q) & \simeq \left \{ \left. d = \pm \prod_{p \in \Sigma} p^{a_p} \in \frac {\mathbb Q^\times}{(\mathbb Q^\times)^2} \, \right \vert \,  \text{$C_d''(\mathbb Q_\nu) \neq \emptyset$ for all $\nu \in \Sigma$} \right \}, \\ \sha(E''/\mathbb Q)[\hat \varphi] & \simeq \biggl \{ \biggl. \{ C_d''/\mathbb Q \} \in WC(E''/\mathbb Q)[\hat \varphi] \, \biggr \vert \,  \text{$C_d''(\mathbb Q_\nu) \neq \emptyset$ for all $\nu \in \Sigma$} \biggr \}. \end{aligned} \end{equation}

\noindent Hence by the canonical exact sequence
\begin{equation} \begin{CD} 1 @>>> \dfrac {(\mathbb Q^\times)^2}{(\mathbb Q^\times)^4} @>>> \dfrac {\mathbb Q^\times}{(\mathbb Q^\times)^4} @>>> \dfrac {\mathbb Q^\times}{(\mathbb Q^\times)^2} @>>> 1 \end{CD} \end{equation}

\noindent we may think of $S^{(\hat \phi)}(E'/\mathbb Q)$ as being a cover of $S^{(\hat \varphi)}(E''/\mathbb Q)$.

\begin{proof} The first statement follows from Proposition \ref{2-isogeny} and the results from \cite{MR87g:11070}, so we focus on the second statement.  Consider the diagram
\begin{equation} \begin{CD} C_{d,t}' @>{\theta}>> E_t' @>{\hat \phi}>> E_t \\ @. @VV{\eta}V @. \\ C_{d,t}'' @>{\vartheta}>> E_t'' @>{\hat \varphi}>> E_t \end{CD} \end{equation} 

\noindent where $\theta$ is the isomorphism in \eqref{theta} and $\vartheta$ is the isomorphism that sends
\begin{equation} \vartheta: \quad (Z,W) \mapsto \left( 2 \, \frac {Z^2 + \sqrt{d} \, W + d}{Z^2}, \ 2 \sqrt{d} \, \frac {(t^2+2) \, Z^2 + 2\sqrt{d} \, W + 2 \, d}{Z^3} \right); \end{equation}

\noindent which has the inverse
\begin{equation} \vartheta^{-1}: \quad (\mathfrak u, \mathfrak v) \mapsto \left( \sqrt{d} \, \frac {2 \, \mathfrak v}{\mathfrak u \, (\mathfrak u - 4)}, \ \sqrt{d} \, \frac {-4 \, t^2 + 2 \, t^2 \, \mathfrak u + {\mathfrak u}^2}{\mathfrak u \, (\mathfrak u - 4)} \right). \end{equation}

\noindent Define the map $\eta_\ast: C_{d,t}' \to C_{d,t}''$ by $\eta_\ast = \vartheta^{-1} \circ \eta \circ \theta$. This gives the composition
\begin{equation} \psi'' \circ \eta_\ast = \left( \hat \varphi \circ \vartheta \right) \circ \left( \vartheta^{-1} \circ \eta \circ \theta \right) = \left( \hat \varphi \circ \eta \right) \circ \theta = \hat \phi \circ \theta = \psi'. \end{equation}

\noindent Moreover, one checks that $\eta_\ast$ maps $(z,w)$ as above.  \end{proof}

\subsection{Modified Descent via 4-Isogeny}

We present an effective version of the four-descent by using the two-descent.  This is the main result of the paper.

\begin{descent} \label{descent} Fix $t \in \mathbb Q^\times$, and denote the elliptic curves
\begin{equation} \begin{aligned} E_t: & \quad & v^2 & = u^3 + \left( t^2 + 2 \right) u^2 + u, \\ E_t': & \quad & V^2 & = U^3 - 2 \left( t^2 - 4 \right) \, U^2 + \left( t^2 + 4 \right)^2 U, \\ E_t'': & \quad  & {\mathfrak v}^2 & = {\mathfrak u}^3 + \left( t^2 - 4 \right) {\mathfrak u}^2 - 4 \, t^2 \, \mathfrak u. \end{aligned} \end{equation}
\begin{enumerate}
\item There are isogenies $\phi: E_t \to E_t'$ of degree 4 with kernel generated by $(-1,t)$, and $\varphi: E_t \to E_t''$ of degree 2 generated by $[2] \, (-1,t) = (0,0)$.
\item Upon identifying the Selmer groups as subgroups of quotients of $\mathbb Q^\times / (\mathbb Q^\times)^4$ and Shafarevich-Tate groups as a collection of homogeneous spaces, there is an exact diagram
\begin{equation} \begin{CD} 0 @>>> \dfrac {E_t(\mathbb Q)}{\hat \phi \left( E_t'(\mathbb Q)\right)} @>{\delta'}>> S^{(\hat \phi)}(E_t' / \mathbb Q) @>>> \sha(E_t' / \mathbb Q)[\hat \phi] @>>> 0 \\ @. @VVV @VVV @VV{\eta}V \\ 0 @>>> \dfrac {E_t(\mathbb Q)}{\hat \varphi \left( E_t''(\mathbb Q)\right)} @>{\delta''}>> S^{(\hat \varphi)}(E_t'' / \mathbb Q) @>>> \sha(E_t'' / \mathbb Q)[\hat \varphi] @>>> 0 \end{CD} \end{equation} 
\noindent where we have the connecting homomorphisms
\begin{equation} \delta': \begin{aligned} (-1,t) & \mapsto - 4 \, t^2 \\ (u,v) & \mapsto u^2 + 2 \, (t^2+1) \, u + 1 - 2 \,t \, v \end{aligned} \qquad \delta'': \quad \begin{aligned} (0,0) & \mapsto  t^2 + 4  \\ (u,v) & \mapsto u \end{aligned} \end{equation}
\noindent while the maps into the Shafarevich-Tate groups send $d \mapsto \{ C_{d,t}' / \mathbb Q \}$ and $d \mapsto \{ C_{d,t}'' / \mathbb Q \}$ in terms of
\begin{equation} \begin{aligned} C_{d,t}': & \quad & d \left( w - \frac {z^2}{4 \, t^2} \right) z^2 = \left( w^2 - d \right)^2, \\ C_{d,t}'': & \quad & d \, W^2 = d^2 + \left( t^2 + 2 \right) d \, Z^2 + Z^4. \end{aligned} \end{equation}
\item The composition
\begin{equation} \begin{CD} \psi': \quad C_{d,t}' @>{\eta}>> C_{d,t}'' @>{\psi''}>> E_t \end{CD} \end{equation} 
\noindent in terms of the maps $\eta$ and $\psi''$ sending
\begin{equation} (z,w) \mapsto \left( -\frac {d \, z^2}{2 \, t \, (w^2 - d)}, \ \frac {d \, z^2 \, (w^2 + d)}{2 \, (w^2 - d)^2} \right), \qquad (Z,W) \mapsto \left( \frac {d}{Z^2}, \ - \frac {d \, W}{Z^3} \right); \end{equation}
\noindent satisfy the congruence $\delta' \bigl(  \psi'(P) \bigr) \equiv d \mod{(\mathbb Q^\times)^4}$.
\end{enumerate}
\end{descent}

Note that we abuse notation and write $\eta$ instead of $\eta_\ast$.

\begin{proof} This follows directly from Theorem \ref{4-descent} and Proposition \ref{2-descent}. \end{proof}


\section{Examples}

In this section we assume that the elliptic curve $E$ and its 4-isogenous curve $E'$ have Weierstrass equations as in Proposition \ref{4-isogeny}.

\subsection{Example: $\mathbb T \simeq \mathbb Z / 4 \mathbb Z$.}

Consider $t = 8$.  Using the ideas in \cite{CG1} as well as Proposition \ref{mordell-weil} it is easy to see that the torsion subgroup of both $E(\mathbb Q)$ and $E'(\mathbb Q)$ is $\mathbb Z / 4 \, \mathbb Z$.  In fact, one may use either Cremona's \texttt{mwrank} or the package \texttt{MAGMA} so see that the Mordell-Weil groups are
\begin{equation} \begin{aligned} E(\mathbb Q) & \simeq \mathbb Z / 4 \mathbb Z \times \mathbb Z, \\ E'(\mathbb Q) & \simeq \mathbb Z / 4 \mathbb Z \times \mathbb Z, \\ E''(\mathbb Q) & \simeq \mathbb Z / 2 \mathbb Z \times \mathbb Z / 2 \mathbb Z \times \mathbb Z. \end{aligned} \end{equation}

\noindent We consider the diagram in \eqref{descent-diagram} in detail.

Performing a 2-descent with \texttt{mwrank} (recall that $\eta: E' \to E''$ and $\hat \varphi: E'' \to E$ are 2-isogenies) one computes the quantities
\begin{equation} \begin{aligned} \dfrac {E''(\mathbb Q)}{\eta \left( E'(\mathbb Q)\right)} & \simeq S^{(\eta)}(E' / \mathbb Q) \simeq \left \{ \pm 1, \, \pm 2 \right \}, & \qquad \sha(E' / \mathbb Q)[\eta] & \simeq \{ 0 \}; \\ \dfrac {E(\mathbb Q)}{\hat \varphi \left( E''(\mathbb Q)\right)} & \simeq S^{(\hat \varphi)}(E'' / \mathbb Q) \simeq \{ \pm 1 \}, & \qquad \sha(E'' / \mathbb Q)[\hat \varphi] & \simeq \{ 0 \}. \end{aligned} \end{equation} 

\noindent Considering the diagram in \eqref{descent-diagram} we see that
\begin{equation} \dfrac {E(\mathbb Q)}{\hat \phi \left( E'(\mathbb Q)\right)} \simeq S^{(\hat \phi)}(E' / \mathbb Q) \simeq \dfrac {\mathbb Z}{2 \, \mathbb Z} \times \dfrac {\mathbb Z}{2 \, \mathbb Z} \qquad \text{and} \qquad \sha(E' / \mathbb Q)[\hat \phi] \simeq \{ 0 \}. \end{equation} 

Using Theorem \ref{4-descent}, we may compute generators explicitly by exploiting the homogeneous spaces.  The curves $E$ and $E'$ both have conductor $2^4 \cdot 17$ so that they have good reduction away from $\Sigma = \{ 2, \, 17, \, \infty \}$.  The homogeneous spaces of interest are of the form
\begin{equation} C_d': \quad d \left( w + a \, z^2 \right) z^2 = \left( w^2 - d \right)^2, \quad a = - 4^{-4}. \end{equation}

\noindent When $d = -1$ we find the solution $(z,w) = (4,0)$; this can be predicted from the connecting homomorphism $\delta': (0,0) \mapsto a^{-1} \equiv -1$.  When $d = 4$ we find the solution $(z,w) = (16, \, 10)$; this maps to the point of infinite order
\begin{equation} \left( \frac w{z^2}, \ \frac {w^2 - d}{z^4} \right) = \left( \frac 5{128}, \, \frac 3{2048} \right) \in E(\mathbb Q). \end{equation}

\noindent Hence $S^{(\hat \phi)}(E'/\mathbb Q) \simeq \{ \pm1, \ \pm 4 \}$.  A brief search shows that $(-1/2, \, -3/4) \in E'(\mathbb Q)$.  In fact, one checks that
\begin{equation} \begin{aligned} \phi \left( 0, \, 0 \right) & = \mathcal O, & \quad \phi \left( \frac 5{128}, \, \frac 3{2048} \right) & = [2] \left( - \frac 12, \ -\frac 34 \right); \\ \hat \phi \left( 0, \, 0 \right) & = [2] \left( 0, \, 0 \right), & \quad \hat \phi \left( -\frac 12, \, -\frac 34 \right) & = [2] \left( \frac 5{128}, \, \frac 3{2048} \right). \end{aligned} \end{equation}

\noindent This shows that
\begin{equation} \frac {E(\mathbb Q)}{\hat \phi \left( E'(\mathbb Q) \right)} \simeq \frac {\mathbb Z}{2 \, \mathbb Z} \times \frac {\mathbb Z}{2 \, \mathbb Z} \qquad \text{and} \qquad \frac {E'(\mathbb Q)}{\phi \left( E(\mathbb Q) \right)} \simeq \frac {\mathbb Z}{2 \, \mathbb Z} \times \frac {\mathbb Z}{4 \, \mathbb Z}. \end{equation}

\subsection{Example: $\mathbb T \simeq \mathbb Z / 2 \mathbb Z \times \mathbb Z / 4 \mathbb Z$.}

Consider $t = 3/2$.  One computes that the Mordell-Weil groups are
\begin{equation} \begin{aligned} E(\mathbb Q) & \simeq \mathbb Z / 2 \mathbb Z \times \mathbb Z / 4 \mathbb Z, \\ E'(\mathbb Q) & \simeq \mathbb Z / 8 \mathbb Z, \\ E''(\mathbb Q) & \simeq \mathbb Z / 2 \mathbb Z \times \mathbb Z / 4 \mathbb Z. \end{aligned} \end{equation}

\noindent We consider Proposition \ref{mordell-weil} in more detail because $\sqrt{ t + \sqrt{t^2 + 4}}$ is rational.  

Performing a 2-descent with \texttt{mwrank} one computes the quantities
\begin{equation} \begin{aligned} \dfrac {E''(\mathbb Q)}{\eta \left( E'(\mathbb Q)\right)} & \simeq S^{(\eta)}(E' / \mathbb Q) \simeq \left \{ \pm 1 \right \}, & \qquad \sha(E' / \mathbb Q)[\eta] & \simeq \{ 0 \}; \\ \dfrac {E(\mathbb Q)}{\hat \varphi \left( E''(\mathbb Q)\right)} & \simeq S^{(\hat \varphi)}(E'' / \mathbb Q) \simeq \{ \pm 1 \}, & \qquad \sha(E'' / \mathbb Q)[\hat \varphi] & \simeq \{ 0 \}; \end{aligned} \end{equation} 

\noindent and so
\begin{equation} \dfrac {E(\mathbb Q)}{\hat \phi \left( E'(\mathbb Q)\right)} \simeq S^{(\hat \phi)}(E' / \mathbb Q) \simeq \dfrac {\mathbb Z}{2 \, \mathbb Z} \qquad \text{and} \qquad \sha(E' / \mathbb Q)[\hat \phi] \simeq \{ 0 \}. \end{equation} 

We now compute explicit points.  The curves have good reduction away from $\Sigma = \{ 2, \, 3, \, 5, \, \infty \}$, and the homogeneous spaces of interest are of the form
\begin{equation} C_d': \quad d \left( w + a \, z^2 \right) z^2 = \left( w^2 - d \right)^2, \quad a = - 3^{-2}. \end{equation}

\noindent When $d = -9$ we find the solution $(z,w) = (3,-3)$; this maps to the point
\begin{equation} \left( \frac w{z^2}, \ \frac {w^2 - d}{z^4} \right) = \left( -\frac 13, \, \frac 29 \right) \in E(\mathbb Q). \end{equation}

\noindent (Note how we may predict $S^{(\hat \phi)}(E'/\mathbb Q) \simeq \{ 1, -9 \}$ from the connecting homomorphism $\delta': (0,0) \mapsto a^{-1} \equiv -9$.)  In fact, the group of rational points $E(\mathbb Q)$ is generated by $\left( -1/3, \, 2/9 \right)$ and $\left(0,0 \right)$, while the group $E'(\mathbb Q)$ is generated by $\left( 25/64, \, 125/1024 \right)$; with orders 2, 4, and 8, respectively.  (Compare with the proof of Proposition \ref{mordell-weil} where $\gamma = \sqrt{ t + \sqrt{t^2 + 4}} = 2$.)  We have $\phi \left( 0, \, 0 \right) = \mathcal O$ while
\begin{equation} \phi \left( - \frac 13, \, \frac 29 \right) = [4] \left( \frac {25}{64}, \ -\frac {125}{1024} \right) \quad \text{and} \quad \hat \phi \left( \frac {25}{64}, \, \frac {125}{1024} \right) =  \left( -\frac 13, \, \frac 29 \right) \oplus (0,0); \end{equation}

\noindent this shows that
\begin{equation} \frac {E(\mathbb Q)}{\hat \phi \left( E'(\mathbb Q) \right)} \simeq \frac {\mathbb Z}{2 \, \mathbb Z} \qquad \text{and} \qquad \frac {E'(\mathbb Q)}{\phi \left( E(\mathbb Q) \right)} \simeq \frac {\mathbb Z}{4 \, \mathbb Z}. \end{equation}

\subsection{Example: $\mathbb T \simeq \mathbb Z / 2 \mathbb Z \times \mathbb Z / 8 \mathbb Z$.}

Consider the collection of elliptic curves studied in \cite{CG1} as discussed in the introduction; each rational $r$ corresponds to $t = \left( r^4 - 6 \, r^2 + 1 \right) / \left( 2 \, r^3 - 2 \, r \right)$.  Bounds on the rank can be computed by \texttt{mwrank} and by \texttt{MAGMA}, but these programs do not seem to be able to compute the rank exactly.  For these five examples we will consider the non-negative integer $R$ such that
\begin{equation} \begin{aligned} E(\mathbb Q) & \simeq \mathbb Z / 2 \mathbb Z \times \mathbb Z / 8 \mathbb Z \times \mathbb Z^R, \\ E'(\mathbb Q) & \simeq \mathbb Z / 4 \mathbb Z \times \mathbb Z^R, \\ E''(\mathbb Q) & \simeq \mathbb Z / 2 \mathbb Z \times \mathbb Z / 4 \mathbb Z \times \mathbb Z^R. \end{aligned} \end{equation}

\noindent Various data about the curves is summarized in Table \ref{curvedata}.  We will see below that, assuming the conjectures of Birch and Swinnerton-Dyer, the lower bounds on $R$ are sharp and the Shafarevich-Tate groups are non-trivial.

\begin{table} \caption{Data for $E_t(\mathbb Q)$ where $t = \left( r^4 - 6 \, r^2 + 1 \right) / \left( 2 \, r^3 - 2 \, r \right)$} \label{curvedata}
\begin{tabular}{|c|cc|c|} 
\hline & & & \\[-11pt] 
$r$ & $t$ & Set of Bad Primes $\Sigma$ & Rank Bounds \\ \hline
& & & \\[-8pt] 
$\dfrac {15}{56}$ & $-\dfrac {5651521}{4890480}$ & $\{ 2, 3, 5, 7, 41, 71, 1231, 3361, 4591, \infty \}$ & $2 \leq R \leq 4$ \\[8pt]
$\dfrac {24}{65}$ & $-\dfrac {3580801}{11384880}$ & $\{ 2, 3, 5, 7, 13, 23, 41, 89, 967, 4801, \infty \}$ & $2 \leq R \leq 4$ \\[8pt]
$\dfrac {11}{69}$ & $-\dfrac {4806319}{1760880}$ & $\{ 2, 3, 5, 7, 11, 23, 29, 223, 2441, 3079, \infty \}$ & $2 \leq R \leq 4$ \\[8pt]
$\dfrac 7{88}$ & $-\dfrac {57695201}{9480240}$ & $\{ 2, 3, 5, 7, 11, 19, 23, 79, 113, 281, 7793, \infty \}$ & $2 \leq R \leq 4$ \\[8pt]
$\dfrac {12}{97}$ & $-\dfrac {80420641}{21568920}$ & $\{ 2, 3, 5, 7, 17, 41, 97, 109, 233, 991, 11593, \infty \}$ & $1 \leq R \leq 4$ \\[8pt] \hline
\end{tabular} \end{table}

We focus on size of the Selmer groups.  The program \texttt{mwrank} computes the groups $S^{(\eta)}(E' / \mathbb Q)$ and $S^{(\hat \varphi)}(E'' / \mathbb Q)$ associated with the 2-isogenies, so it suffices to consider the group $S^{(\hat \phi)}(E' / \mathbb Q)$ associated with the 4-isogeny.  We may compute the size of the latter Selmer group using the exact diagram in \eqref{descent-diagram}; the data is collected in Table \ref{groups}.  Consider the image of the torsion subgroup of $E(\mathbb Q)$ under the connecting homomorphism $\delta': (x,y) \mapsto x^2 - y \mod{(\mathbb Q^\times)^4}$.  Quite explicitly, the torsion subgroup is generated by the points
\begin{equation} \begin{aligned} P & = \left( \frac {r^2}{( r^2 + 2 \, r - 1) \, (r^2 - 2 \, r - 1)}, \ \frac {2 \, r^4}{( r^2 + 2 \, r - 1)^2 \, (r^2 - 2 \, r - 1)^2} \right) \\ Q & = \left( \frac {r \, (r+1) \, (r-1)^2}{( r^2 + 2 \, r - 1) \, (r^2 - 2 \, r - 1)^2}, \ \frac {r \, (r+1)^2 \, (r-1)^3}{( r^2 + 2 \, r - 1)^2 \, (r^2 - 2 \, r - 1)^3} \right) \end{aligned} \end{equation}

\noindent having orders 2 and 8, respectively, and the connecting homomorphism $\delta'$ sends
\begin{equation} P \mapsto - \frac {r^4}{( r^2 + 2 \, r - 1)^2 \, (r^2 - 2 \, r - 1)^2}, \qquad Q \mapsto \frac {r \, (r+1)^3 \, (r-1)^3}{( r^2 + 2 \, r - 1)^2 \, (r^2 - 2 \, r - 1)^4}. \end{equation} 

\noindent For the $d \in S^{(\hat \phi)}(E' / \mathbb Q)$ above, consider the homogeneous spaces
\begin{equation} C_d': \quad d \left( w + a \, z^2 \right) z^2 = \left( w^2 - d \right)^2, \quad a = - \frac {r^2 \, (r+1)^2 \, (r-1)^2}{( r^2 + 2 \, r - 1)^2 \, (r^2 - 2 \, r - 1)^2}. \end{equation}

\noindent Some generators $(x,y)$ on $E(\mathbb Q)$ can be found by \texttt{mwrank}; these can be pulled back to a point $(z,w)$ on $C_d'(\mathbb Q)$ via the mapping $\psi': C_d' \to E$ in Theorem \ref{4-descent}.  We summarize these rational points in Table \ref{points}.

\begin{table} \caption{Selmer Group Computations via \texttt{mwrank}} \label{groups}
\begin{tabular}{|c|cc|c|} 
\hline & & & \\[-11pt] 
$r$ & $S^{(\eta)}(E' / \mathbb Q)$ & $S^{(\hat \varphi)}(E'' / \mathbb Q)$ & $\# S^{(\hat \phi)}(E' / \mathbb Q)$ \\[3pt] \hline
& & & \\[-8pt] 
$\dfrac {15}{56}$& $\langle -1, 2, 3, 5, 7, 41, 71 \rangle$ & $\langle -1, 2, 3, 5, 7, 41, 71 \rangle$ & $2^{13}$ \\[8pt]
$\dfrac {24}{65}$ & $\langle -1, 2, 3, 5, 13, 41, 89 \rangle$ & $\langle -1, 2, 3, 5, 13, 41, 89 \rangle$ & $2^{13}$ \\[8pt]
$\dfrac {11}{69}$ & $\langle -1, 2, 5, 23, 29, 33 \rangle$ & $\langle -1, 2, 5, 23, 29, 33 \rangle$ & $2^{11}$ \\[8pt]
$\dfrac 7{88}$ & $\langle -1, 2, 7, 11, 15, 57, 113, 843 \rangle$ & $\langle -1, 2, 7, 57, 11 \cdot 15 \rangle$ & $2^{12}$ \\[8pt]
$\dfrac {12}{97}$ & $\langle -1, 2, 51, 291, 1635, 57965 \rangle$ & $\langle -1, 2, 51, 291, 1635 \rangle$ & $2^{10}$ \\[8pt] \hline
\end{tabular} \end{table}

\begin{table} \caption{Points $(z,w) \in C_d'(\mathbb Q)$ Corresponding to $d \in S^{(\hat \phi)}(E'/\mathbb Q)$} \label{points}
\begin{tabular}{|c|ccc|} 
\hline & & & \\[-11pt] 
$r$ & $d$ & $z$ & $w$ \\[3pt] \hline
& & & \\[-8pt] 
$\dfrac {15}{56}$ & $-1 \cdot 5651521^2$ & $\dfrac {5651521}{840}$ & $5651521$ \\[8pt]
& $2 \cdot 8733^2$ & $\dfrac {86033006559349680248647}{367497397145772723536}$ & $\dfrac {3354273204372689787801}{320119683924889132}$ \\[8pt]
& $41 \cdot 71 \cdot 8617^2$ & $\dfrac {32245866415137}{37085066480}$ & $\dfrac {1001789642807}{3730892}$ \\[8pt]
& $3\cdot 5 \cdot 7 \cdot 56515210^2$ & $\dfrac {37902870933931733579}{3188191909720998}$ & $\dfrac {5724171826291316638240}{12960129714313}$ \\[8pt] \hline
& & & \\[-8pt] 
$\dfrac {24}{65}$ & $-1 \cdot 6769^2$ & $\dfrac {155687}{1560}$ & $6769$ \\[8pt] 
& $2 \cdot 15^2$ & $\dfrac {3893880282213551}{2376079357686320}$ & $\dfrac {45660470553}{2504457868}$ \\[8pt] 
& $41 \cdot 89 \cdot 130^2$ & $\dfrac {83760182508961}{72984831998502}$ & $\dfrac {26347553866400}{3333554033}$ \\[8pt]
& $3 \cdot 5 \cdot 13 \cdot 527982^2$ & $\dfrac {4690464211822333}{4980114669368}$ & $\dfrac {128551710083251455}{15183276431}$ \\[8pt] \hline
& & & \\[-8pt] 
$\dfrac {11}{69}$ & $-1 \cdot 4806319^2$ & $\dfrac {9612638}{759}$ & $4806319$ \\[8pt] 
& $5 \cdot 4683^2$ & $\dfrac {68687104829}{107909560}$ & $\dfrac {585885447}{42652}$ \\[8pt] 
& $23 \cdot 2453^2$ & $\dfrac {670894550748541}{5454390952560}$ & $\dfrac {71568676933}{12390148}$ \\[8pt]
& $2 \cdot 29 \cdot 33 \cdot 215530^2$ & $\dfrac {1001805034739772531298429}{834126406602556183512}$ & $\dfrac {24018346018297942811880275}{3159569721979379483}$ \\[8pt] \hline
& & & \\[-8pt] 
$\dfrac 7{88}$ & $-1 \cdot 57695201^2$ & $\dfrac {57695201}{616}$ & $57695201$ \\[8pt]
& $2 \cdot 6441^2$ & $\dfrac {5596434497}{2776356}$ & $\dfrac {41962889}{738}$ \\[8pt] 
& $7 \cdot 7658655^2$ & $\dfrac {15608602011772334}{143349694803}$ & $\dfrac {271628289087723675}{2275391981}$ \\[8pt]
& $11 \cdot 15 \cdot 57 \cdot 135723^2$ & $\dfrac {26233163006176556}{5388839579079}$ & $\dfrac {5056981351335585}{954960053}$ \\[8pt] \hline
& & & \\[-8pt] 
$\dfrac {12}{97}$ & $-1 \cdot 80420641^2$ & $\dfrac {80420641}{1164}$ & $80420641$ \\[8pt] 
& $17 \cdot 97 \cdot 34^2$ & $\dfrac {3173285824541959}{59077245130310}$ & $\dfrac {1304349940207}{558755747}$ \\[8pt]
& $1635 \cdot 34779^2$ & $\dfrac {115031523259990561}{28757077855714}$ & $\dfrac {31422940575120}{79995877}$ \\[8pt] \hline
\end{tabular} \end{table}

Rathbun \cite{rathbun:01} has remarked that the points in Table \ref{points} correspond to the complete list generators on the elliptic curve if we assume the validity of the Birch and Swinnerton-Dyer conjectures.  That is, we consider the order of the vanishing of the $L$-series in order to compute the rank, as well as the relationship with the residue and regulator of the elliptic curve to compute the generators; the information is summarized in Table \ref{bsd}.  (For more information on the BSD conjectures, consult \cite[Chapter 4]{MR2002g:11073}.)  Notice that the $\hat \phi$-component of the Shafarevich-Tate group is not always of square order; we give an explanation.  For the moment, let $\phi: E \to E'$ be \emph{any} isogeny.  The canonical alternating, non-degenerate, bilinear pairing
\begin{equation} \begin{CD} \sha(E' / \mathbb Q)  \times \sha(E' / \mathbb Q) @>>> \mathbb Q / \mathbb Z \end{CD} \end{equation}

\begin{table} \caption{Shafarevich-Tate Group Computations Assuming BSD} \label{bsd}
\begin{tabular}{|c|c|c|c|} 
\hline & & & \\[-11pt] 
$r$ & Rank & $\# S^{(\hat \phi)}(E' / \mathbb Q)$ & $\# \sha(E' / \mathbb Q)[\hat \phi]$ \\[3pt] \hline
& & & \\[-8pt] 
$\dfrac {15}{56}$& 2 & $2^{13}$ & $2^6$ \\[8pt]
$\dfrac {24}{65}$ & 2 & $2^{13}$ & $2^6$ \\[8pt]
$\dfrac {11}{69}$ & 2 & $2^{11}$ & $2^4$ \\[8pt]
$\dfrac 7{88}$ & 2 & $2^{12}$ & $2^5$ \\[8pt]
$\dfrac {12}{97}$ & 1 & $2^{10}$ & $2^7$ \\[8pt] \hline
\end{tabular} \end{table}

\noindent induces a non-degenerate, bilinear pairing
\begin{equation} \begin{CD} \sha(E' / \mathbb Q)[\hat \phi]  \times \dfrac {\sha(E' / \mathbb Q)}{\phi \left( \sha(E / \mathbb Q) \right)} @>>> \mathbb Q / \mathbb Z. \end{CD} \end{equation}

\noindent (For more information on the Cassels pairing, consult \cite[Appendix C, \S 17]{MR87g:11070}.)  When $\phi$ is the ``multiplication-by-$m$'' map we have $\sha(E / \mathbb Q) / m \, \sha(E / \mathbb Q) \simeq \sha(E / \mathbb Q)[m]$ -- assuming that $\sha(E / \mathbb Q)$ is finite -- so the pairing induces an alternating, non-degenerate, bilinear pairing
\begin{equation} \begin{CD} \sha(E / \mathbb Q)[m] \times \sha(E / \mathbb Q)[m] @>>> \mathbb Q / \mathbb Z. \end{CD} \end{equation}

\noindent Hence $\sha(E / \mathbb Q)[m]$ must have square order.  In general, however, $\hat \phi$ might not be self-dual; that is, upon considering the short exact sequence
\begin{equation} \begin{CD} 0 @>>> \dfrac {\sha(E / \mathbb Q)}{\sha(E / \mathbb Q)[\phi]} @>{\phi}>> \sha(E'/\mathbb Q) @>>> \dfrac {\sha(E' / \mathbb Q)}{\phi \left( \sha(E / \mathbb Q) \right)} @>>> 0 \end{CD} \end{equation}

\noindent we might have
\begin{equation} \# \frac {\sha(E / \mathbb Q)}{\sha(E / \mathbb Q)[\phi]} \neq \# \frac {\sha(E' / \mathbb Q)}{\sha(E' / \mathbb Q)[\hat \phi]} \implies \frac {\sha(E' / \mathbb Q)}{\phi \left( \sha(E / \mathbb Q) \right)} \not \simeq \sha(E' / \mathbb Q)[\hat \phi]. \end{equation}


\bibliographystyle{plain}

\end{document}